\documentclass[12pt]{article}

\usepackage{amsfonts,graphicx,euscript}

\newtheorem{theorem}{Theorem}
\newtheorem{corollary}{Corollary}

\newcommand{\trace}{\mathrm{trace}\:}
\newcommand{\conv}{\mathrm{conv}\:}

\title{Semidefinite representation of\\
convex hulls of rational varieties}

\author{Didier Henrion$^{1,2}$}
\date{\today}

\begin{document}

\maketitle

\footnotetext[1]{CNRS; LAAS; 7 avenue du colonel Roche, F-31077 Toulouse; France;
Universit\'e de Toulouse; UPS, INSA, INP, ISAE; LAAS; F-31077 Toulouse; France}
\footnotetext[2]{Faculty of Electrical Engineering, Czech Technical University in Prague,
Technick\'a 2, CZ-16626 Prague, Czech Republic}
\addtocounter{footnote}{2}

\begin{abstract}
Using elementary duality properties of positive semidefinite moment matrices
and polynomial sum-of-squares decompositions, we prove that the
convex hull of rationally parameterized algebraic varieties is semidefinite
representable (that is, it can be represented as a projection of an affine
section of the cone of positive semidefinite matrices) in the case of (a) curves;
(b) hypersurfaces parameterized by quadratics; and (c) 
hypersurfaces parameterized by bivariate quartics;
all in an ambient space of arbitrary dimension.
\end{abstract}

\section{Introduction}

Semidefinite programming, a versatile extension of linear programming
to the convex cone of positive semidefinite matrices (semidefinite
cone for short), has found many
applications in various areas of applied mathematics and engineering,
especially in combinatorial optimization, structural mechanics
and systems control. For example, semidefinite programming was
used in \cite{hsk} to derive linear matrix inequality (LMI)
convex inner approximations of non-convex semi-algebraic
stability regions, and in \cite{hl} to derive a hierarchy of embedded
convex LMI outer approximations of non-convex semi-algebraic sets
arising in control problems.

It is easy to prove that affine sections and projections of the
semidefinite cone are convex semi-algebraic sets, but it is still
unknown whether all convex semi-algebraic sets can be modeled
like this, or in other words, whether all convex semi-algebraic sets
are semidefinite representable. Following the development of
polynomial-time
interior-point algorithms to solve semidefinite programs, a long list
of semidefinite representable semi-algebraic sets and convex hulls
was initiated in \cite{nn} and completed in \cite{bn}. Latest
achievements in the field are reported in \cite{l} and \cite{hn}.

In this paper we aim at enlarging the class of semi-algebraic sets
whose convex hulls are explicitly semidefinite representable.
Using elementary duality properties of positive semidefinite moment
matrices and polynomial sum-of-squares decompositions -- nicely
recently surveyed in \cite{laurent} -- we prove that the
convex hull of rationally parameterized algebraic varieties is
explicitly semidefinite representable in the case of (a) curves;
(b) hypersurfaces parameterized by quadratics; and (c) 
hypersurfaces parameterized by bivariate quartics;
all in an ambient space of arbitrary dimension.

Rationally parameterized surfaces arise often in engineering,
and especially in computer-aided design (CAD). For example,
the CATIA (Computer Aided Three-dimensional Interactive Application)
software, developed since 1981 by the French company Dassault Syst\`emes,
uses rationally parameterized surfaces as its core 3D surface
representation. CATIA was originally used to develop Dassault's
Mirage fighter jet for the French airforce, and then it was adopted
in aerospace, automotive, shipbuilding, and other industries.
For example, Airbus aircrafts are designed in Toulouse with the
help of CATIA, and architect Frank Gehry has used the software to
design his curvilinear buildings, like the Guggenheim Museum
in Bilbao or the Dancing House in Prague, near the Charles Square
buildings of the Czech Technical University.

\section{Notations and definitions}

Let $x = [x_0,\:x_1,\:,\cdots,x_m] \in {\mathbb R}^{m+1}$ and
\[
\zeta_d(x) = [x_0^d,\:x_0^{d-1}x_1,\: x_0^{d-1}x_2,\:\cdots,x_0^{d-2}x_1^2,\:\cdots,x_m^d]
\in {\mathbb R}^{s(m,d)}[x]
\]
denote a basis vector of $m$-variate forms of degree $d$, with
$s(m,d) = (m+d)!/(m!d!)$.
Let $y = [y_{\alpha}]_{|\alpha|\leq 2d} \in {\mathbb R}^{s(m,2d)}$
be a real-valued sequence indexed in basis $\zeta_{2d}(x)$,
with $\alpha \in {\mathbb N}^m$ and $|\alpha|=\sum_k \alpha_k$.
A form $x \mapsto p(x)
= p^T \zeta_{2d}(x)$ is expressed in this basis via its coefficient
vector $p \in {\mathbb R}^{s(m,2d)}$. Given a sequence $y \in {\mathbb R}^{s(m,2d)}$,
define the linear mapping $p \mapsto L_y(p) = p^T y$, and
the linear moment matrix $M_d(y)$ satisfying the relation $L_y(pq) = p^TM_d(y)q$ for all
$p,q \in {\mathbb R}^{s(m,d)}$. It has entries
$[M_d(y)]_{\alpha,\beta} = L_y([\zeta_d(x)\zeta_d(x)^T]_{\alpha,\beta})
= y_{\alpha+\beta}$ for all $\alpha,\beta \in {\mathbb N}^m$,
$|\alpha|+|\beta|\leq 2d$.
For example, when $m=2$ and $d=2$ (trivariate quartics) we have
$s(m,2d)=15$. To the form $p(x)=x_0^4-x_0x_1x_2^2+5x_1^3x_2$ we
associate the linear mapping $L_y(p)=y_{00}-y_{12}+5y_{31}$.
The 6-by-6 moment matrix is given by
\[
M_2(y) = \left[\begin{array}{c|cc|ccc}
y_{00} & * & * & * & * & * \\ \hline
y_{10} & y_{20} & * & * & * & * \\
y_{01} & y_{11} & y_{02} & * & * & * \\ \hline
y_{20} & y_{30} & y_{21} & y_{40} & * & * \\
y_{11} & y_{21} & y_{12} & y_{31} & y_{22} & * \\
y_{02} & y_{12} & y_{03} & y_{22} & y_{13} & y_{04}
\end{array}\right]
\]
where symmetric entries are denoted by stars.
See \cite{laurent} for more details on these notations
and constructions.

Given a set $\mathcal Z$, let $\conv{\mathcal Z}$
denote its convex hull, the smallest convex
set containing $\mathcal Z$. Finally, the notation
$M_d(y) \succeq 0$ means that matrix $M_d(y)$ is
positive semidefinite.

\section{Convex cones and moment matrices}

Consider the Veronese variety
\[
{\mathcal W}_{m,d} = \{\zeta_{2d}(x) \in {\mathbb R}^{s(m,2d)}
\: :\: x \in {\mathbb R}^{m+1}\}
\]
and the convex cones
\[
{\mathcal Z}_{m,d} = {\conv}{\mathcal W}_{m,d} 
\]
and
\[
{\mathcal Y}_{m,d} =
\{y \in {\mathbb R}^{s(m,2d)} \: : \: M_d(y) \succeq 0\}.
\]

\begin{theorem}\label{zeta}
If $m=1$ or $d=1$ or $d=m=2$
then ${\mathcal Z}_{m,d} = {\mathcal Y}_{m,d}$.
\end{theorem}

{\bf Proof}:
The inclusion ${\mathcal Z}_{m,d} \subset {\mathcal Y}_{m,d}$
follows from the definition of a moment matrix
since
\[
M_d(\zeta_{2d}(x)) = \zeta_d(x)\zeta_d(x)^T \succeq 0.
\]
The converse inclusion is shown by contradiction.
Assume that $y^* \notin {\mathcal Z}_{m,d}$ and hence that
there exists a (strictly separating) hyperplane $\{y \: : p(y)=0\}$
such that $p^T y^* < 0$ and $p^T y \geq 0$ for all $y \in {\mathcal Z}_{m,d}$.
It follows that form $x \mapsto p(x) = p^T\zeta_{2d}(x)$ is
globally non-negative. Since $m=1$ or $d=1$ or $d=m=2$, the
form can be expressed as a sum of squares of forms
\cite[Theorem 3.4]{laurent} and we can write $p(x) = \sum_k q_k^2(x)
= \sum_k (q_k^T \zeta_d(x))^2 
= \zeta_d(x)^T P \zeta_d(x)$ for some matrix $P = \sum_k q_k q_k^T
\succeq 0$. Then $L_y(p) = p^T y = \trace(PM_d(y)) = \sum_k q_k^T M_d(y) q_k$.
Since $L_{y^*}(p) < 0$, there must be an index $k$ such that 
$q_k^T M_d(y^*) q_k < 0$ and hence matrix $M_d(y^*)$ cannot
be positive semidefinite, which proves that $y^* \notin {\mathcal Y}_{m,d}$.
$\Box$

See also \cite{fialkow} for a study of the moment problem
in the bivariate quartic case ($d=m=2$).

\section{Rational varieties}

Given a matrix $A \in {\mathbb R}^{(n+1)\times s(m,2d)}$, we define
the rational variety ${\mathcal V}_{m,d}$ (of degree $2d$ with $m$ parameters in an $n$-dimensional
ambient space) as an affine projection of the Veronese variety ${\mathcal W}_{m,d}$:
\[
{\mathcal V}_{m,d} = {\mathcal A}({\mathcal W}_{m,d}) = 
\{v \in {\mathbb R}^n \: :\: \left[\begin{array}{c}1\\v\end{array}\right] = A \zeta_{2d}(x), \:
x \in {\mathbb R}^{m+1}\}.
\]
Theorem \ref{zeta} identifies the cases when
the convex hull of this rational variety is exactly semidefinite
representable. That is, when it can be formulated as
the projection of an affine section of the semidefinite cone.

\begin{corollary}\label{exact}
If $m=1$ or $d=1$ or $d=m=2$ then 
\[
\conv{\mathcal V}_{m,d} = \{v \in {\mathbb R}^n \: :\:
\left[\begin{array}{c}1\\v\end{array}\right] =  Ay, \:
M_d(y) \succeq 0, \: y \in {\mathbb R}^{s(m,2d)}\}.
\]
\end{corollary}

{\bf Proof:}
We have $\conv {\mathcal V}_{m,d} = \conv {\mathcal A}({\mathcal W}_{m,d})
= {\mathcal A}(\conv{\mathcal W}_{m,d}) = {\mathcal A}({\mathcal Z}_{m,d})$
and the result follows readily from Theorem \ref{zeta}.$\Box$

The case $m=1$ corresponds to rational curves.
The case $d=1$ corresponds to quadratically parameterized rational hypersurfaces.
The case $d=m=2$ corresponds to hypersurfaces parameterized by bivariate
quartics. All these rational varieties live in an ambient space of arbitrary
dimension $n>m$.

In all other cases, the inclusion $\conv{\mathcal V}_{m,d} \subset
{\mathcal A}({\mathcal Y}_{m,d})$ is strict.
For example, when $d=3,\:m=2$, the vector $y^* \in {\mathbb R}^{28}$
with non-zero entries
\[
y_{00}^*=32,\:y_{20}^*=y_{02}^*=34,
\:y_{40}^*=y_{04}^*=43,\:y_{22}^*=30,\:
y_{60}^*=y_{06}^*=128,\:y_{42}^*=y_{24}^*=28
\]
is such that $M_3(y^*)\succ 0$
but $L_y^*(p^*)<0$
for the Motzkin form $p^*(x) = x_0^6-3x_0^2x_1^2x_2^2+x_1^4x_2^2+x_1^2x_2^4$
which is globally non-negative.
In other words, $y^* \in {\mathcal A}({\mathcal Y}_{m,d})$
but $y^* \notin \conv{\mathcal V}_{m,d}$.

\section{Examples}

\subsection{Parabola}

The parabola
\[
{\mathcal V} = \{v \in {\mathbb R}^2 \: :\:
v_1^2 - v_2 = 0\}
\]
can be modeled as an affine projection
of a quadratic Veronese variety
\[
{\mathcal V} = \{v \in {\mathbb R}^2 \: :\:
\left[\begin{array}{c}1\\v_1\\v_2\end{array}\right] =
\left[\begin{array}{c}x_0^2\\x_0x_1\\x_1^2\end{array}\right], \:
x \in {\mathbb R}^2\},
\]
i. e. $n=2$, $d=1$, $m=1$ and $A$ is the 3-by-3 identity matrix
in the notations of the previous section.

By Corollary \ref{exact}, the convex hull of the
parabola is the set
\[
\begin{array}{rcl}
\conv{\mathcal V} & = & \{v \in {\mathbb R}^2 \: :\:
\left[\begin{array}{c}1\\v_1\\v_2\end{array}\right] =
\left[\begin{array}{c}y_0\\y_1\\y_2\end{array}\right], \:
M_1(y) = \left[\begin{array}{cc}y_0&y_1\\y_1&y_2\end{array}\right] \succeq 0, \:
y \in {\mathbb R}^3\} \\ \\
& = & \{v \in {\mathbb R}^2 \: :\:
\left[\begin{array}{cc}1&v_1\\v_1&v_2\end{array}\right] \succeq 0\}
\end{array}
\]
which is described with a 2x2 LMI.

\subsection{Trefoil knot}

\begin{figure}[h!]
\begin{center}
\includegraphics[width=8cm]{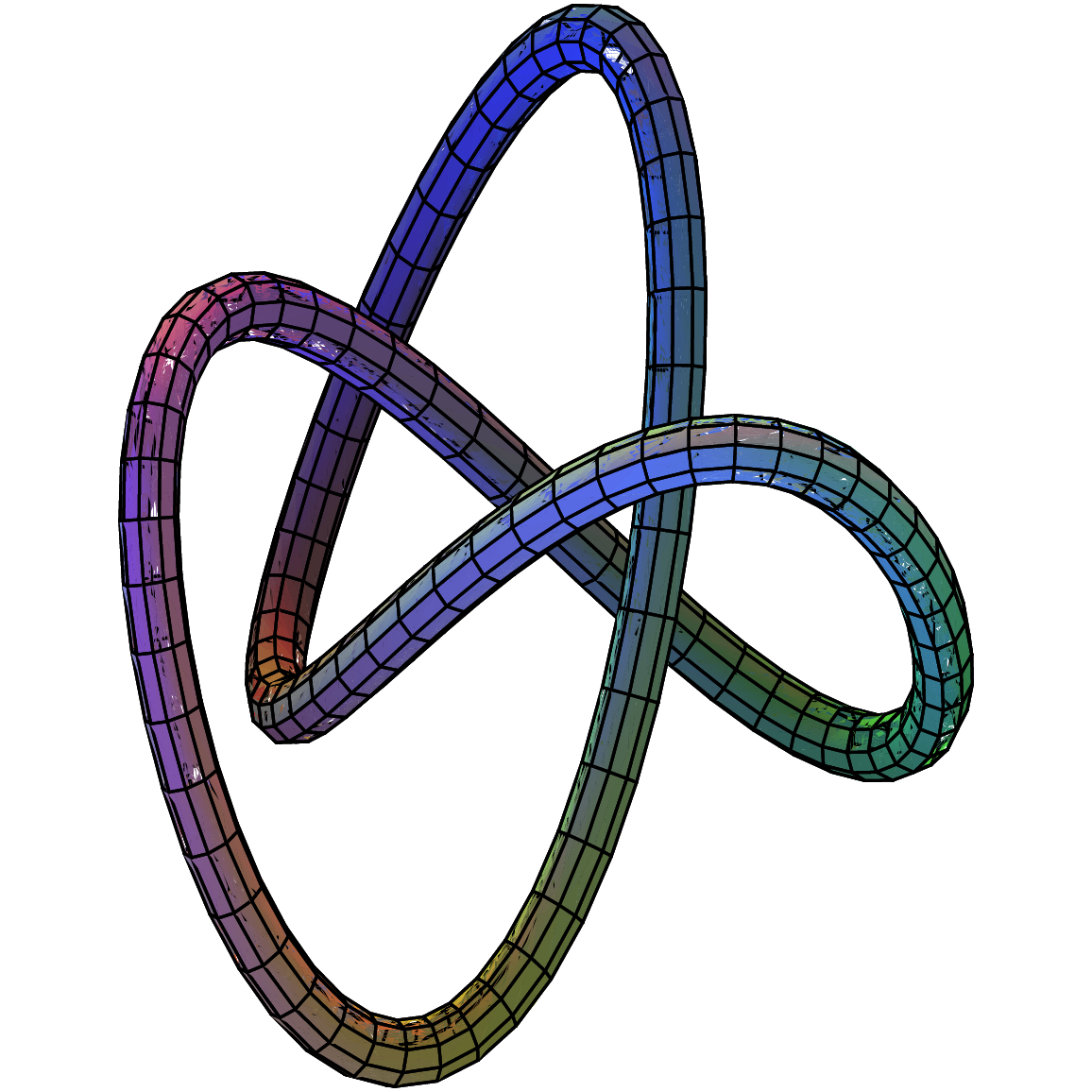}
\caption{Tube plot of the trefoil knot curve, whose
convex hull is exactly semidefinite representable
with 3 liftings.\label{trefoilknot}}
\end{center}
\end{figure}

The space trigonometric curve
\[
{\mathcal V} = \{v \in {\mathbb R}^3 \: :\:
v_1(\alpha) = \cos\alpha+2\cos2\alpha, \:
v_2(\alpha) = \sin\alpha+2\sin2\alpha, \:
v_3(\alpha) = 2\sin3\alpha, \:
\alpha \in [0,2\pi]\}
\]
is called a trefoil knot, see \cite{bk} and
Figure \ref{trefoilknot}.

Using the standard change of variables
\[
\cos\alpha = \frac{x_0^2-x_1^2}{x_0^2+x_1^2}, \quad
\sin\alpha = \frac{2x_0x_1}{x_0^2+x_1^2}
\]
and trigonometric formulas,
the space curve admits a rational representation
as an affine projection of a sextic Veronese variety
\[
\begin{array}{ll}
{\mathcal V} = \{v \in {\mathbb R}^3 \: :\: &
1 = (x_0^2+x_1^2)^3, \:
v_1 = (x_0^2+x_1^2)(3x_0^4-12x_0^2x_1^2+x_1^4), \:\\ &
v_2 = 2x_0x_1(x_0^2+x_1^2)(5x_0^2-3x_1^2), \:
v_3 = 4x_0x_1(x_0^2-3x_1^2)(3x_0^2-x_1^2), \:\:
x \in {\mathbb R}^2\}
\end{array}
\]
i.e. $n=3$, $m=1$ and $d=3$ in the notations
of the previous section.

By Corollary \ref{exact}, the convex hull of the
trefoil knot curve is exactly semidefinite representable as
\[
\conv{\mathcal V} = \{v \in {\mathbb R}^3 \: :\:
\left[\begin{array}{c}1\\v\end{array}\right] = A y, \:
M_3(y) \succeq 0, \: y \in {\mathbb R}^7\}
\]
with 
\[
A = \left[\begin{array}{rrrrrrr}
1 & 0 & 3 & 0 & 3 & 0 & 1 \\
3 & 0 & -9 & 0 & -11 & 0 & 1 \\
0 & 10 & 0 & 4 & 0 & -6 & 0 \\
0 & 12 & 0 & -40 & 0 & 12 & 0
\end{array}\right]
\]
and
\[
M_3(y) = \left[\begin{array}{cccc}
y_0 & * & * & * \\
y_1 & y_2 & * & * \\
y_2 & y_3 & y_4 & * \\
y_3 & y_4 & y_5 & y_6
\end{array}\right]
\]
where symmetric entries are denoted by stars.
The affine system of equations involving $v$ and $y$
can be solved by Gaussian elimination to yield
the equivalent formulation:
\[
\begin{array}{l}
\conv{\mathcal V} = \{v \in {\mathbb R}^3 \: :\: \\
\quad
\left[\begin{array}{cccc}
\frac{1}{6}(3+v_1+2u_1-4u_3) & * & * & *\\
\frac{1}{112}(10v_2+v_3+48u_2) &
\frac{1}{18}(3-v_1-20u_1-2u_3) & * & * \\
\frac{1}{18}(3-v_1-20u_1-2u_3) &
\frac{1}{224}(6v_2-5v_3+96u_2) &
u_1 & * \\
\frac{1}{224}(6v_2-5v_3+96u_2) &
u_1 &
u_2 & 
u_3
\end{array}\right] \succeq 0, 
\:\: u \in {\mathbb R}^3\}
\end{array}
\]
which is an explicit semidefinite representation with
3 liftings.

\subsection{Steiner's Roman surface}

\begin{figure}[h!]
\begin{center}
\includegraphics[width=7cm]{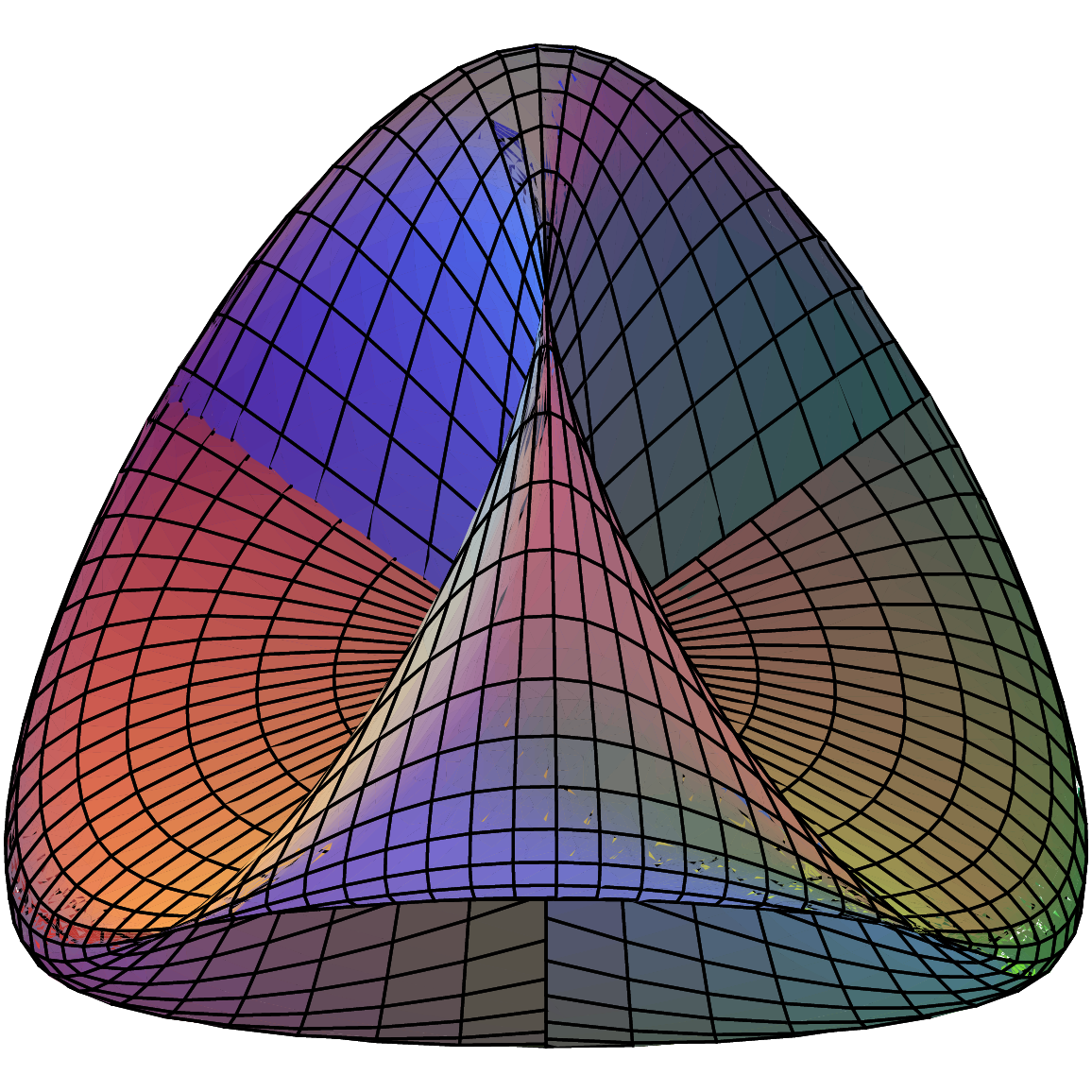} 
\hspace{5mm}
\includegraphics[width=7cm]{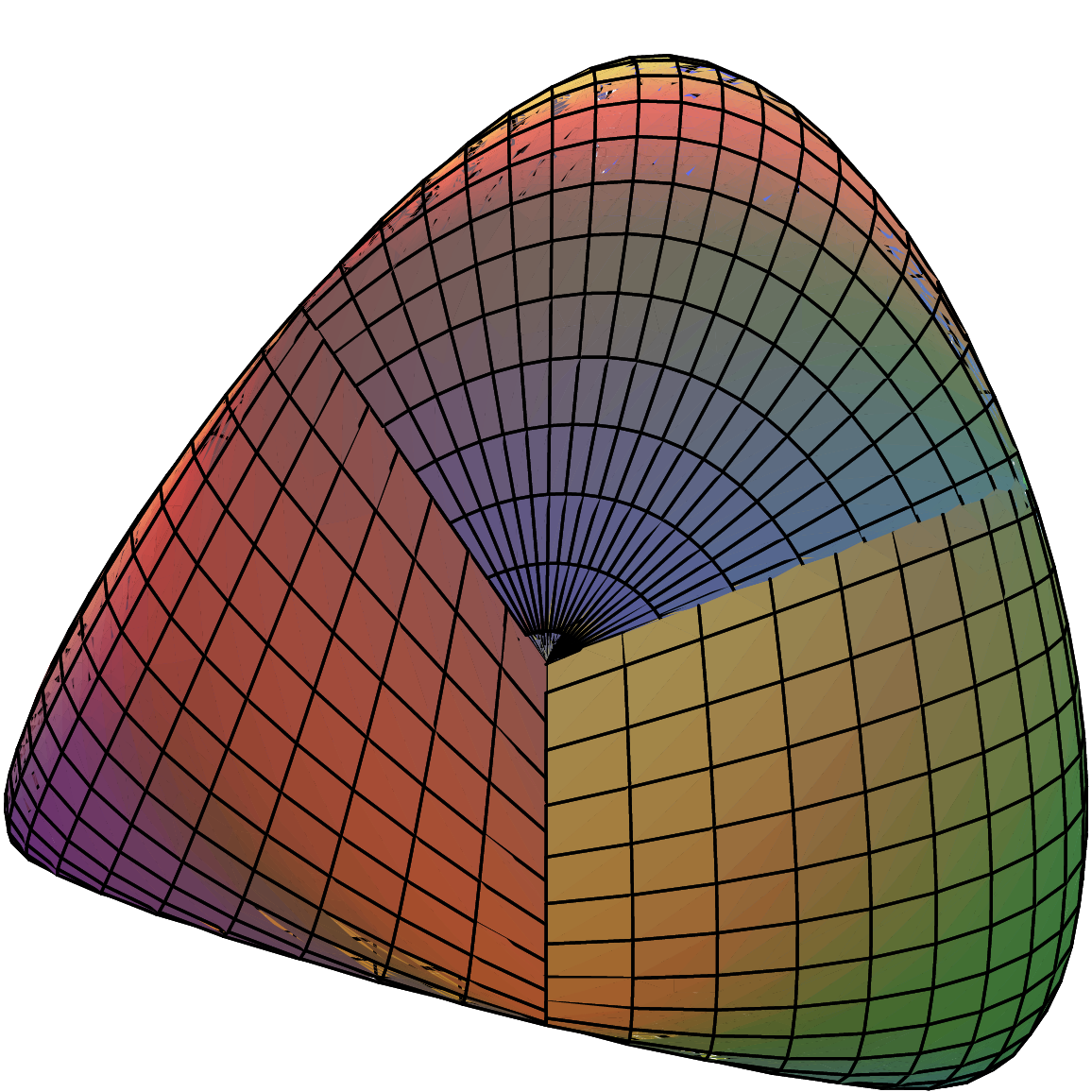}
\caption{Two views of Steiner's Roman surface, whose
convex hull is semidefinite representable with 2 liftings.
\label{roman}}
\end{center}
\end{figure}

Quadratically parameterizable rational surfaces are classified in \cite{csc}.
A well-known example is Steiner's Roman surface, a non-orientable
quartic surface with three double lines, which is 
parameterized as follows:
\[
{\mathcal V} = \{v \in {\mathbb R}^3 \: :\:
v_1 = \frac{2x_1}{1+x_1^2+x_2^2}, \:
v_2 = \frac{2x_2}{1+x_1^2+x_2^2}, \:
v_3 = \frac{2x_1x_2}{1+x_1^2+x_2^2}, \:\:
x \in {\mathbb R}^2\}
\]
see Figure \ref{roman}.

The surface can be modeled as an affine projection
of a quadratic Veronese variety:
\[
{\mathcal V} = \{v \in {\mathbb R}^3 \: :\:
1 = x_0^2+x_1^2+x_2^2, \:
v_1 = 2x_0x_1, \:
v_2 = 2x_0x_2, \:
v_3 = 2x_1x_2, \:\:
x \in {\mathbb R}^3\}
\]
i.e. $n=3$, $m=2$ and $d=1$ in the notations
of the previous section. By Corollary \ref{exact}, its convex hull is exactly
semidefinite representable as
\[
\conv{\mathcal V} = \{v \in {\mathbb R}^3 \: :\:
\left[\begin{array}{c}1\\v\end{array}\right] = A y, \:
M_1(y) \succeq 0, \: y \in {\mathbb R}^6\}
\]
with 
\[
A = \left[\begin{array}{rrrrrr}
1 & 0 & 0 & 1 & 0 & 1 \\
0 & 2 & 0 & 0 & 0 & 0 \\
0 & 0 & 2 & 0 & 0 & 0 \\
0 & 0 & 0 & 0 & 2 & 0
\end{array}\right]
\]
and
\[
M_1(y) = \left[\begin{array}{ccc}
y_{00} & * & * \\
y_{10} & y_{20} & * \\
y_{01} & y_{11} & y_{02} \\
\end{array}\right].
\]
The affine system of equations involving $v$ and $y$ can
easily be solved to yield
the equivalent formulation:
\[
\conv{\mathcal V} = \{v \in {\mathbb R}^3 \: :\: 
\left[\begin{array}{ccc}
1-u_1-u_2 & * & * \\
\frac{1}{2}v_1 & u_1 & * \\
\frac{1}{2}v_2 & \frac{1}{2}v_3 & u_2
\end{array}\right] \succeq 0, \:\:
u \in {\mathbb R}^2\}
\]
which is an explicit semidefinite representation with
2 liftings.

\subsection{Cayley cubic surface}

\begin{figure}[h!]
\begin{center}
\includegraphics[width=8cm]{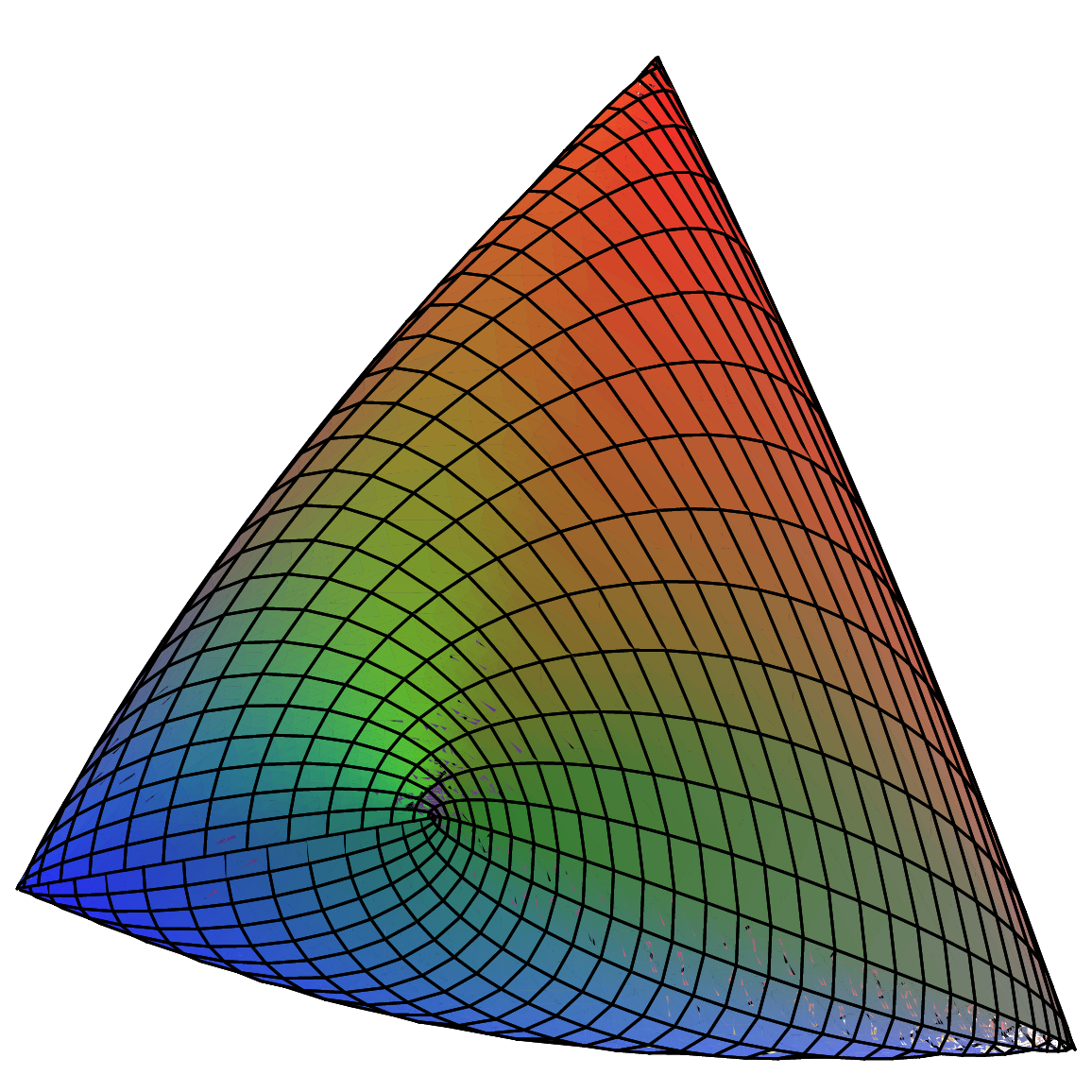}
\caption{Convex connected component of Cayley's cubic surface,
semidefinite representable with no liftings.\label{cayley}}
\end{center}
\end{figure}

Steiner's Roman surface, studied in the
previous paragraph, is dual to Cayley's cubic surface
$\{v \in {\mathbb R}^3 \: :\: \det C(v) = 0\}$
where
\[
C(v) = \left[\begin{array}{ccc}
1 & * & * \\
v_1 & 1 & * \\
v_2 & v_3 & 1
\end{array}\right].
\]
The origin belongs to a set delimited by a convex
connected component of this surface, admitting the following
affine trigonometric parameterization:
\[
\begin{array}{ll}
{\mathcal V} = \{v \in {\mathbb R}^3 \: : &
v_1(\alpha) = \cos\alpha_1, \:
v_2(\alpha) = \sin\alpha_2, \\
& v_3(\alpha) = \cos\alpha_1 \sin\alpha_2 - \cos\alpha_2 \sin\alpha_1, \:
\alpha_1 \in [0,\pi], \:
\alpha_2 \in [-\pi,\pi]\}.
\end{array}
\]
This is the boundary of the LMI region
\[
\conv{\mathcal V} = \{v \in {\mathbb R}^3 \: :\: C(v) \succeq 0\}
\]
which is therefore semidefinite representable with no liftings.
This set is a smoothened
tetrahedron with four singular points, see Figure \ref{cayley}.

Using the standard change of variables
\[
\cos\alpha_i = \frac{x_0^2-x_i^2}{x_0^2+x_i^2}, \quad
\sin\alpha_i = \frac{2x_0x_i}{x_0^2+x_i^2}, \quad i=1,2
\]
we obtain an equivalent rational parameterization
\[
\begin{array}{ll}
{\mathcal V} = \{v \in {\mathbb R}^3 \: :\: &
1 = (x_0^2+x_1^2)(x_0^2+x_2^2), \:
v_1 =(x_0^2-x_1^2)(x_0^2+x_2^2), \\
& v_2 = 2x_0x_2(x_0^2+x_1^2), \:
v_3 = 2x_0(-x_1+x_2)(x_0^2+x_1x_2), \:\:
x \in {\mathbb R}^3\}.
\end{array}
\]
which is an affine projection of a quadratic
Veronese variety, i.e.
$n=3$, $m=2$ and $d=2$ in the notations
of the previous section.
By Corollary \ref{exact}, its convex hull is exactly
semidefinite representable as
\[
\conv{\mathcal V} = \{v \in {\mathbb R}^3 \: :\:
\left[\begin{array}{c}1\\v\end{array}\right] = A y, \:
M_2(y) \succeq 0, \: y \in {\mathbb R}^{15}\}
\]
with $A$ of size $4$-by-$15$
and $M_2(y)$ of size $6$-by-$6$,
not displayed here. It follows that $\conv{\mathcal V}$
is semidefinite representable as a $6$-by-$6$ LMI
with $11$ liftings.

We have seen however that $\conv{\mathcal V}$
is also semidefinite representable as a $3$-by-$3$
LMI with no liftings, a considerable simplification.
It would be interesting to design an algorithm simplifying
a given semidefinite representation, lowering the size
of the matrix and the number of variables. As far as
we know, no such algorithm exists at this date.

\section{Conclusion}

The well-known equivalence between polynomial non-negativity
and existence of a sum-of-squares decomposition was used,
jointly with semidefinite programming duality, to
identify the cases for which the convex hull of a rationally
parameterized variety is exactly semidefinite representable.
Practically speaking, this means that optimization of a linear
function over such varieties is equivalent to semidefinite
programming, at the price of introducing a certain number
of lifting variables.

If the problem of detecting whether a plane algebraic curve
is rationally parameterizable, and finding explicitly such
a parametrization, is reasonably well understood
from the theoretical and numerical point of view -- see
\cite{swp} and M. Van Hoeij's {\tt algcurves} Maple package
for an implementation -- the case of surfaces is
much more difficult \cite{schicho}. Up to our knowledge,
there is currently no working computer implementation of a
parametrization algorithm for surfaces. Since an explicit
parametrization is required for an explicit semidefinite
representation of the convex hull of varieties,
the general case of algebraic varieties given in implicit
form (i.e. as a polynomial equation), remains largely open.

Finally, we expect that these semidefinite representability
results may have applications when studying non-convex
semi-algebraic sets and varieties arising from stability
conditions in systems control, in the spirit of \cite{hsk,hl}.
These developments are however out of the scope of the present paper.

\section*{Acknowledgments}

The first draft benefited from technical advice by Jean-Bernard Lasserre,
Monique Laurent, Josef Schicho and two anonymous reviewers.
I am grateful to Bernd Sturmfels for pointing out an error in the coefficients
of my original lifted LMI formulation of the trefoil curve
when preparing paper \cite{sturmfels}.

\end{document}